\documentclass[3p,times,twocolumn]{elsarticle} 

\usepackage{float}
\usepackage{graphicx}
\usepackage{url}
\usepackage{amssymb}
\usepackage[svgnames]{xcolor}
\usepackage{array}
\usepackage{booktabs}
\usepackage{colortbl}
\usepackage{multirow}
\usepackage{tabularx}
\usepackage{amsfonts}
\usepackage{cancel}
\usepackage{empheq}
\usepackage{bm}
\usepackage{parskip}
\usepackage{framed}
\usepackage{mathtools}
\usepackage{units}
\usepackage{refcount}

\usepackage[bottom]{footmisc}
\usepackage[hidelinks]{hyperref}

\usepackage{cleveref}
\crefname{equation}{}{}
\Crefname{equation}{}{}

\crefname{theo}{Theorem}{Theorems}
\Crefname{theo}{Theorem}{Theorems}
\crefname{prop}{Proposition}{Propositions}
\Crefname{prop}{Proposition}{Propositions}
\crefname{assm}{Assumption}{Assumptions}
\Crefname{assm}{Assumption}{Assumptions}
\crefname{lemm}{Lemma}{Lemmas}
\Crefname{lemm}{Lemma}{Lemmas}
\crefname{defn}{Definition}{Definitions}
\Crefname{defn}{Definition}{Definitions}
\crefname{remk}{Remark}{Remarks}
\Crefname{remk}{Remark}{Remarks}
\crefname{cond}{Condition}{Conditions}
\Crefname{cond}{Condition}{Conditions}

\usepackage{autonum}

\usepackage{enumitem}
\setlist[itemize]{leftmargin=*,noitemsep,align=left}
\setlist[enumerate]{leftmargin=*,noitemsep,align=left}

\usepackage{algorithm}
\usepackage{algorithmicx}
\usepackage{algpseudocode}
\algrenewcommand\algorithmicindent{3mm}
\makeatletter

\usepackage{amsthm}
\makeatletter
\def\thm@space@setup{%
  \thm@preskip=\parskip \thm@postskip=0pt
}
\makeatother

\newcommand*{\preccur}{\includegraphics[scale=0.4]{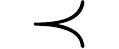}}

\theoremstyle{plain}
\newtheorem{theo}{Theorem}
\newtheorem{prop}[theo]{Proposition}
\newtheorem{assm}[theo]{Assumption}
\newtheorem{lemm}[theo]{Lemma}

\theoremstyle{definition}
\newtheorem{defn}[theo]{Definition}
\newtheorem{remk}[theo]{Remark}
\newtheorem*{remk*}{Remark}
\newtheorem{cond}[]{Condition}

\newtheorem{innercustomthm}{Condition}
\newenvironment{ccond}[1]
  {\renewcommand\theinnercustomthm{#1}\innercustomthm}
  {\endinnercustomthm}

\journal{Control Engineering Practice}

\begin{document}

\begin{frontmatter}

\title{Analytical results for the multi-objective design of model-predictive control}

\author[1]{Vincent Bachtiar\corref{cor1}}\ead{bachtiarv@unimelb.edu.au}
\author[1]{Chris Manzie}\ead{manziec@unimelb.edu.au}
\author[1]{William H. Moase}\ead{moasew@unimelb.edu.au}
\author[2,3]{Eric C. Kerrigan}\ead{e.kerrigan@imperial.ac.uk}

\address[1]{Dept. of Mechanical Engineering, The University of Melbourne, VIC 3010, Australia}
\address[2]{Dept. of Electrical and Electronic Engineering, Imperial College London, London SW7 2AZ, U.K.}
\address[3]{Dept. of Aeronautics, Imperial College London, London SW7 2AZ, U.K.}

\cortext[cor1]{Corresponding author}

\begin{abstract}
In model-predictive control (MPC), achieving the best closed-loop performance under a given computational resource is the underlying design consideration. This paper analyzes the MPC design problem with control performance and required computational resource as competing design objectives. The proposed multi-objective design of MPC (MOD-MPC) approach extends current methods that treat control performance and the computational resource separately -- often with the latter as a fixed constraint -- which requires the implementation hardware to be known a priori. The proposed approach focuses on the tuning of structural MPC parameters, namely sampling time and prediction horizon length, to produce a set of optimal choices available to the practitioner. The posed design problem is then analyzed to reveal key properties, including smoothness of the design objectives and parameter bounds, and establish certain validated guarantees. Founded on these properties, necessary and sufficient conditions for an effective and efficient solver are presented, leading to a specialized multi-objective optimizer for the MOD-MPC being proposed. Finally, two real-world control problems are used to illustrate the results of the design approach and importance of the developed conditions for an effective solver of the MOD-MPC problem.
\end{abstract}


\end{frontmatter}

\section{Introduction} \label{sc01}

Model-predictive control (MPC) is a typically computationally expensive method of approaching the control of constrained systems. As a result, the computational resource required at each sampling instant is a consideration in the overall design process. This is particularly true in systems with fast dynamics, where there is often significant conflict between the complexity of the problem considered at each time step and the available time to find a solution. The close interrelation between control performance and required computational resource warrants that these indices are analyzed in synchrony to streamline the design process and avoid unnecessary costs. Both objectives depend on a number of tuning parameters of the optimal control problem including, but not limited to, the sampling time, prediction horizon length, and fidelity/order of the prediction model.


Previously, much focus has been given to find the best control performance in a single-objective optimization design problem, separate to the consideration of the required computational resource. However, there are still a number of knowledge gaps in existing MPC design approaches. MPC tuning for control performance is mostly done via methods that rely on rules-of-thumb and general guidelines \cite{Garriga:10,Qin:03,Rani:97}. Further developments have been made consequently, based on metaheuristics such as particle swarm optimization \cite{Nery:14} and genetic algorithms \cite{Lee:08}, as well as gradient descent \cite{Bunin:12}, for the single-objective optimization of MPC.

Several multi-objective optimization approaches for control system design have also been studied for the optimization of control performance. Similar to that of the single-objective counterpart, metaheuristic methods are prevalently used for the multi-objective tuning of classical control, such as PID \cite{Ayala:12,Reynoso:13,Xue:10}, sliding mode control \cite{Mahmoodabadi:14,Taherkhorsandi:14}, as well as others \cite{Reynoso:14}. A similar approach is applied in MPC tuning by using an off-the-shelf method of goal attainment \cite{Exadaktylos:10,Vega:08}. Although more systematic than general guidelines, these methods provide non-specialized approaches that do not exploit certain characteristics of the problem and potentially require a rather exhaustive and possibly computationally impractical search to produce an optimal design set. As an alternative to approaches based on guidelines and metaheuristics, analytical methods employing problem simplifications have been proposed \cite{Bagheri:14,Shah:11,Shridhar:98}. However, these typically overlook some aspects of the original problem such as explicit constraint handling.


The studies discussed so far consider control performance as the sole design objective, whether with a single- or multi-objective outlook. The approach separates software and hardware design, revealing only half the insight in control design. Hardware design largely determines the implementation cost of the controller and is often not known a priori, thus is a part of the design process. The co-design of software and hardware provides a more comprehensive approach that optimizes control performance, as well as implementation cost that is dictated by the required computational resource to functionally implement the control system. Rather than treating the required resource as a fixed constraint, it should be co-optimized alongside control performance, avoiding system over-design or the need to re-design the system. Furthermore, previous studies (e.g.\ \cite{Bagheri:14,Exadaktylos:10,Nery:14}) have typically assumed the structural parameters of the controller -- such as sampling rate and prediction horizon -- are fixed. Nonetheless, structural MPC parameters have been shown to have an underlying role for MPC design improvement \cite{Bachtiar:15,Bachtiar:16}.

In light of the above discussion, the value of a co-design approach in streamlining the design process of control systems has been noted \cite{Allison:14}. Further, the fundamental concept of a software and hardware co-design approach for real-time optimization has been studied \cite{Kerrigan:14}, although analytical results to support applications in MPC are still yet to be fully developed. The main contribution of this paper is a systematic development of the optimal MPC design with a multi-objective approach. Theoretical results concerning the nature of the design problem are presented to establish certain assumptions and guarantees. These results are then used to understand the nature of the optimization problem at hand and subsequently provide conditions that a selected optimizer must satisfy in order to effectively and efficiently compute the optimal (Pareto) frontier. The approach allows the practitioner to understand the trade-off between performance and resources in structurally designing an MPC controller for a given real-world control problem.


The paper is outlined as follows; Section~\ref{sc02} contains the MPC formulation studied. The proposed multi-objective MPC design approach is then presented in Section~\ref{sc03}. Section~\ref{sc04} identifies the key properties of the multi-objective problem, including smoothness properties and parameter bounds. In Section~\ref{sc05}, conditions for an effective and efficient solver are presented and a compliant algorithm is proposed. Section~\ref{sc06} considers two real-world examples to demonstrate the design approach and importance of the conditions developed for an effective solver. Section~\ref{sc07} presents conclusions of the study and potential future work.

\subsection*{Notational conventions and definitions}


$\|\mathbf{v}\|^2_M\coloneqq \mathbf{v}^\mathsf{T}M \mathbf{v}$. $\otimes$ and $\oslash$ denotes element-wise multiplication and division, respectively. $\mathsf{U}[a,b]$ is a random number uniformly distributed between $a$ and $b$. Unless stated otherwise, an ordered list (column vector) is defined with a bold typeface e.g.\ $\mathbf{v}$ with its size denoted by $|\mathbf{v}|$. The element values are thus $\mathbf{v}\coloneqq\left(v_1,\ldots,v_{|\mathbf{v}|}\right)$. A set containing several ordered list is defined in calligraphy e.g.\ $\mathcal{V}$ with entries $\mathcal{V}\coloneqq\left\{\mathbf{v}_1,\ldots,\mathbf{v}_{|\mathcal{V}|}\right\}$.

\section{Controller design} \label{sc02}

Consider a nonlinear dynamic plant model
\begin{equation}
	\dot{\mathsf{x}} = \mathsf{f}(\mathsf{x},\mathsf{u})
\end{equation}
with states $\mathsf{x}(t)\in\mathbb{R}^{n_\mathsf{x}}$ and inputs $\mathsf{u}(t)\in\mathbb{R}^{n_\mathsf{u}}$ which satisfy standard properties as described in the following.
\begin{assm}
	$(\mathsf{x},\mathsf{u})\mapsto\mathsf{f}(\mathsf{x},\mathsf{u})$ is continuous in $(\mathsf{x},\mathsf{u})$ and globally Lipschitz continuous in $\mathsf{x}$ uniformly in $\mathsf{u}$.
\label[assm]{th07}\end{assm}
\begin{assm}
	$(\mathsf{x},\mathsf{u})\mapsto\mathsf{f}(\mathsf{x},\mathsf{u})$ is differentiable with respect to $\mathsf{u}$ for all $\mathsf{x}\in\mathbb{R}^{n_\mathsf{x}}$.
\label[assm]{th08}\end{assm}

Discretization is used for the purpose of digital control, such that the plant is controlled in a sampled-data fashion at sampling instants $t_i\coloneqq ih$ for $i\in\mathbb{N}_{\geq0}$ with sampling period $h$. The control command sequence is restricted to a zero-order-hold
\begin{equation}
	\mathsf{u}(t) = \mathsf{u}_i,\; \forall t\in[ih,ih+h), i\in\mathbb{N}_{\geq0}.
\end{equation}
The aim is to control the plant by applying a control law $\kappa$ to regulate the model to the origin. The control law depends on the current state $\mathsf{x}_i\coloneqq \mathsf{x}(t_i)$ and the control design parameters $\mathbf{p}$,
\begin{equation}
	\mathsf{u}_i = \kappa(\mathsf{x}_i,\mathbf{p}).
\end{equation}
Let $\mathbf{p}\coloneqq(p_1,\ldots,p_{n_p})$	 contain the design parameters $p_1,\ldots,p_{n_p}$ to be tuned.


In this paper, the control command is obtained by solving a finite-horizon, optimal control problem (OCP) at each sampling instant $t_i$,
\begin{subequations}\noeqref{eq01a,eq01b,eq01c,eq01d,eq01e}
\begin{flalign}
	\hspace{1mm} (x^*(\cdot),u^*(\cdot))\coloneqq\arg\underset{(x,u)}{\min}\; J(x,u,\mathbf{p}) &&
\label{eq01a}
\end{flalign}\vspace{-8mm}
\begin{align}
	\text{s.t. } 
	x(0) 			 		&= \mathsf{x}_i \label{eq01b}\\
	\dot{x}(\tau)	&= Ax + Bu && \forall \tau \in [0,T]\label{eq01c}\\
	x(\tau) 				&\in [\underline{x},\overline{x}],\; u(\tau) \in [\underline{u},\overline{u}] && \forall \tau\in[0,T)\label{eq01d}\\
	u(\tau) 				&= u(kh),\;\forall k \in \mathbb{N}_{\geq0} && \forall \tau\in[kh,kh+h). \label{eq01e}
\end{align}
\label{eq01}\end{subequations}
For succinctness, the dependence of $\tau\mapsto x^*(\tau)$ and $\tau\mapsto u^*(\tau)$ on $(\mathsf{x}_i,\mathbf{p})$ is omitted. Consequently,
\begin{equation}
	\kappa(\mathsf{x}_i,\mathbf{p}) \coloneqq u^*(0). \label{eq17}
\end{equation}
The real-time variable $\mathsf{x}$ is distinct from the predicted variable $x$ used internally in the OCP, although sized equally such that $x(t)\in\mathbb{R}^{n_\mathsf{x}}$ and inputs $u(t)\in\mathbb{R}^{n_\mathsf{u}}$. Further, also note the distinction between the true plant $\mathsf{f}$ and prediction model in linear time-invariant (LTI) form $f \coloneqq Ax + Bu \coloneqq\left.\frac{\partial \mathsf{f}}{\partial\mathsf{x}}\right|_{0,0}x + \left.\frac{\partial \mathsf{f}}{\partial\mathsf{u}}\right|_{0,0}u $ used internally in the OCP. The two models have the same equilibrium at the origin, that is $f(0,0) = \mathsf{f}(0,0)$. The optimization is subject to the prediction model \eqref{eq01c} representing the dynamics of the plant initialized at \eqref{eq01b}, and the plant constraints \eqref{eq01d}. The zero-order-hold control \eqref{eq01e} discretizes the control command over the sampling steps $k\in\{0,\ldots,N-1\}$.

One common choice for the OCP cost function in \eqref{eq01a} is a quadratic
\begin{equation}
	J(x,u,\mathbf{p}) \coloneqq \int_0^{Nh} \|x(\tau)\|^2_Q + \|u(\tau)\|^2_R \;d\tau + \|x(T)\|_{Q_\text{f}}^2
	\label{eq04}
\end{equation}
to penalize the state/input deviations from zero. This cost is composed by the \textit{stage cost} weighted by $Q\geq0$ and $R>0$, and the \textit{terminal cost} weighted by $Q_\text{f}\geq0$. $Nh =: T$ is the prediction horizon length associated with $N$ prediction steps. Finally, the OCP is assumed to be non-degenerate.
\begin{assm}[Non-degeneracy]\label[assm]{th09}
	The OCP \eqref{eq01} is non-degenerate so that its solution $(u^*,x^*)$ is unique.
\end{assm}

\subsection{Design parameters} \label{sc02.1}

From the OCP formulation, a key design parameter is the sampling time $h$ that dictates how often a new control input can be commanded to the plant. This also sets an upper-bound on the time available for the computing hardware to solve the OCP. Next, along with the sampling time, the number of prediction steps governs the length of the prediction horizon $T \coloneqq Nh$. This is the time horizon in which constraints can be applied in the prediction of the future plant behavior, and thus its value affects the performance of the controller. Further, the number of prediction steps directly affects the size of the OCP problem, that is the number of unknowns in the problem. Each of $h$ and $T$ affect both design objectives.

Generally, the prediction model type is a design parameter that can be chosen, for example, as a linear-time invariant (LTI), linear time-varying, or nonlinear model. In this study, the prediction model \cref{eq01c} is defined to be LTI, a particularly application-relevant choice that reduces the general nonlinear OCP into a quadratic program (QP) for which many practical solvers exist. Furthermore, this study focuses on the MPC (software) parameters, so that those that exclusively are attributes of the hardware, such as data precision \cite{Kerrigan:14}, are not considered even if they affect both objectives.

Cost function attributes, namely the cost weighting matrices, are also design parameters. These affect the OCP solution and hence the control performance, and also the time taken to numerically solve the OCP particularly when certain solver algorithms are used. For solvers which are insensitive to ill-conditioning (e.g.\ interior point method \cite{Dang:16}), the parameters can be assumed to only affect control performance. On the other hand, the algorithm used to solve the OCP is a design choice that only affects computational complexity. As long as the solver is convergent, it is assumed that it will find the one local (thus global) minimum of the OCP as a convex problem (QP). A related parameter is the representation of the OCP, e.g.\ dense and sparse representations, for which some are suited to a particular algorithm and some, another. This also can be assumed to produce the same OCP solution and not affect control performance. Finally, a numerical tolerance can be used as an algorithm attribute, dictating the accuracy of the numerical solution of the OCP. This tolerance value would affect both control performance and required computational resource.

In multi-objective design, focus shall be given to design parameters that are \textit{coupled} i.e.\ those that affect both design objectives. In the proposed design problem, the coupled parameters of the MPC architecture, namely the sampling time $h$ and number of prediction steps $N$, are the considered design parameters. Conversely, parameters that only affect one of the objectives, hence \textit{decoupled}, are fixed. For the fixed parameters, in particular the solver algorithm, tolerance and numerical precision, it is assumed that they are well-chosen so that the true, global solution of the OCP can be obtained.
\begin{assm}\label[assm]{th22}
	The numerical solution of the OCP obtained is close to the true/analytical solution.
\end{assm}

 

\subsection{Design objectives} \label{sc02.2}

\subsubsection*{Control performance}

The control performance measures how well the controller steer the states to the origin over time, indicating, for example, the speed of response to set point changes, multi-variable decoupling and damping performances. This can be obtained through a simulated environment of the controlled plant. There are many measures that can be used to indicate control performance, the most common being the integrated squared error (ISE) of the states and inputs,
\begin{subequations}
\begin{align}\hphantom{\ref{eq02a}\ref{eq02b}\ref{eq02c}\ref{eq02d}}\hspace{-16mm}
	U \left(\mathsf{x}_0,\mathbf{p}\right) &\coloneqq \int^\infty_0 v(\mathsf{x}(\tau),\mathsf{u}(\tau)) \;d\tau \label{eq02a}\\
	\text{s.t. }	
	\mathsf{x}(0) 			&= \mathsf{x}_0 \label{eq02b}\\
	\dot{\mathsf{x}} 	&= \mathsf{f}(\mathsf{x},\mathsf{u})\label{eq02c}\\
	\mathsf{u}(\tau) 	&= \kappa^*(\mathsf{x}(t_i),\mathbf{p}), \, \forall \tau\in[ih,ih+h), \, \forall i \in \mathbb{N}_{\geq0}\label{eq02d}
\end{align}\label{eq02}
\end{subequations}
\hspace{-2mm} where the cost function is chosen as
\begin{equation}
	v(\mathsf{x},\mathsf{u}) \coloneqq \|\mathsf{x}(\tau)\|^2_Q + \|\mathsf{u}(\tau)\|^2_R
\end{equation}
in accordance with the OCP cost function \eqref{eq04}. $U\left(\mathsf{x}_0,\mathbf{p}\right)$ defines the closed-loop value function for a given initial condition and design parameter choice.

The feasibility and stability of the closed-loop system shall be introduced as follow.
\begin{assm}[Recursive feasibility]\label[assm]{th24}
	The closed-loop system \eqref{eq02} is recursively feasible such that for a given initial state $\mathsf{x}_0\in \mathcal{X}_\text{\normalfont S}$, the OCP \eqref{eq01} is feasible and remains feasible at all subsequent sampling steps. $\mathcal{X}_\text{\normalfont S}$ is the associated feasible set.
\end{assm}

\vspace{-0.6mm}

\begin{assm}[Stability]\label[assm]{th23}
	The OCP \eqref{eq01} guarantees the asymptotic stability of the closed-loop system \eqref{eq02}.
\end{assm}
The assumptions allow the simulated closed-loop system to run indefinitely and will not encounter an abrupt termination (e.g. due to infeasibility of the OCP). Consequently, the ISE measure \cref{eq02a} is well-defined.

\begin{remk}\label[remk]{th26}
Not all choices of the parameters $(h,N)$ would validate \cref{th23,th24}. The assumptions serve as a starting guideline when choosing the parameters for the closed-loop system to be functional (feasible and stable). The set containing such parameters is to be narrowed further to define those that are optimal.
\end{remk}

Specifying the exact conditions to guarantee \cref{th23,th24} is not the focus of this study and the reader shall refer to the relevant references that are amply available. A brief guideline is given below.
\vspace{-2mm}
\begin{enumerate}
	\item \cite{Levis:71} notes that there exist critical sampling periods where the discrete plant loses full controllability, which consequently leads to an unstable closed-loop trajectory. Since fast systems are of interest, short sampling periods can be appropriately considered (shorter than the first critical sampling periods).
	\item For $f=\mathsf{f}$, an appropriate choice of $Q_\text{f}$ \cite{Mayne:14} and $T$ \cite{Kerrigan:00} is able to guarantee recursive feasibility and stability. The restriction on model fidelity might be relaxed so that the discrepancy between the two is bounded, instead of zero. 
\end{enumerate}
\vspace{-2mm}
If at either one of \cref{th23,th24} is not valid, a penalised value of the metric $U$ can be set for cases when the simulation terminates early. 

A linear combination of the closed-loop value function
\begin{align}
	V(\mathcal{X}_0,\mathbf{p}) &\coloneqq \sum^{|\mathcal{X}_0|}_{i=1} w_i U(\mathsf{x}_{0,i},\mathbf{p}) &
	\forall \mathsf{x}_{0}&\in \mathcal{X}_0
\label{eq10}
\end{align}
provides a numerical measure of the control performance of the MPC law \eqref{eq17} with the OCP \eqref{eq01}. A linear combination is chosen to preserve the smoothness properties of the value-function $U$. $\mathcal{X}_0$ is a set of initial conditions representative of the intended operating range of the controlled plant, each of which is weighted by $w_i$ to determine the relative significance of each scenario.

\subsubsection*{Required computational resource}

In designing the controller, the capability of the computational hardware onto which the controller will be implemented is often not known a priori. It is therefore imperative to set the required computational resource, which is a primary factor in determining cost, as a design objective. Structural MPC parameters dictate the computational complexity of the OCP \eqref{eq01}. In turn, the complexity governs the required resource, affecting the time taken by a processing unit to generate a control command that is upper-bounded by the sampling period.

In many cases, a numerical simulation is conducted instead of directly testing the control plant to test the performance of the control architecture. The computational data obtained in simulation is reflective of the hardware used in the simulation platform, namely the \textit{simulation hardware}. This is to be differentiated from the \textit{implementation hardware}, which is the actual hardware used to implement the controller and control the plant.

Let the upper-bound for the time taken for a command input to be generated using the simulation hardware be $\gamma \left(\mathbf{p}\right)$. This indicates the computational complexity for a given design $\mathbf{p}$. Based on this, a dimensionless measure denoted the \textit{Resource Number} can be derived as
\begin{equation}
	\eta(\mathbf{p}) \coloneqq \gamma(\mathbf{p})/h
\end{equation}
to indicate the required power of the implementation hardware relative to the simulation hardware so that the MPC controller can be functionally implemented.


For the Resource Number to be a meaningful measure, an assumption which specifies the relationship between the simulation and implementation hardware is needed.
\begin{assm}[Scalability]\label[assm]{th17}
	The solution time upper-bound of the simulation hardware $\eta$ is linearly scalable to the solution-time upper-bound $\eta_\text{\normalfont I}$ if the implementation hardware was used. That is, $\eta = a\eta_\text{\normalfont I}$ for some constant multiplier $a\in\mathbb{R}_{>0}$.
\end{assm}
The importance of the scalability assumption \cref{th17} will be clarified later (\cref{th19}) when the multi-objective design problem is formulated.

The solution time upper-bound $\gamma$ is modeled to take a polynomial form as given in the following.
\begin{assm}[Solution time upper-bound]\label[assm]{th20}
The upper-bound on solution time is monotonically increasing with the number of prediction steps $N$, modeled by a polynomial of degree $n$
\begin{equation}
	\gamma(\mathbf{p}) \coloneqq \sum_{i=0}^n a_iN^i
\label{eq16}\end{equation}
for some constants $a_i$, $i\in\{0,\ldots,n\}$ depending on the OCP representation and solver used. 
\end{assm}

The assumption is based on the fact that as $N$ increases, the number of unknowns in the OCP increases as well, which is assumed to further extend to the time needed to solve the OCP.

Generally, the parameters $a_i$ are dependent on design parameters such as sampling time, OCP representation and solver algorithm. For instance, sampling time has an effect on the matrices of the OCP and how well they are conditioned, thus affecting solution time. However, variations in solution time from difference in sampling time are not significant in the tests conducted in this study and are treated as negligible. Furthermore, since other design parameters are treated as constants, the parameters $a_i$ are constant.

\cref{th20} is later verified in the results in Section~\ref{sc06}. Consequently,
\begin{equation}
	\eta(\mathbf{p})\coloneqq \frac{\gamma(\mathbf{p})}{h} \coloneqq \frac{1}{h}\sum_{i=0}^n a_iN^i.
\label{eq09}
\end{equation}

\section{Multi-objective design of MPC} \label{sc03}

Let the objectives be defined as a vector $\boldsymbol\ell\coloneqq(V,\eta)\in\mathbb{R}^{n_\ell}$ where $n_\ell=2$. The multi-objective design (MOD) of MPC is posed as the following:
\begin{subequations}
\begin{align}
	\mathcal{P}_\bullet (\mathcal{P}_\text{s}) \coloneqq \; & \arg\underset{\mathbf{p}}{\text{m-min}}\;\boldsymbol\ell(\mathbf{p}) \label{eq03a} {\color{white}\eqref{eq03a}\eqref{eq03b}}\\
	&\;\text{s.t. }\;\mathbf{p}\in\mathcal{P}_\text{s}. \label{eq03b}
\end{align}\label{eq03}
\end{subequations}
\hspace{-1mm}The minimization (denoted m-min) is a multi-objective minimization to find the Pareto optimal design set $\mathcal{P}_\bullet$ for a given search space $\mathcal{P}_\text{s}$ (see \cref{th26}). This solution set contains the Pareto optimal design choices for the design engineer to select from, based on the Pareto front
\begin{equation}
	\mathcal{L}(\mathcal{P}_\text{s}) \coloneqq \big\{\boldsymbol\ell(\mathbf{p})\,|\,\mathbf{p}\in\mathcal{P}_\bullet(\mathcal{P}_\text{s})\big\}
\label{eq07}\end{equation}
that shows the optimal trade-off between the two objectives. The Pareto optimal (non-dominated) points, defined below, make up the Pareto front.

\begin{defn}[Pareto optimal point \cite{Deb:01}]\label[defn]{th15}
A point $\boldsymbol\ell(\mathbf{p}_\bullet)$ with $\mathbf{p}_\bullet\in \mathcal{P}$ is a Pareto point iff there does not exist another design choice $\mathbf{p}\in \mathcal{P}$ such that $\boldsymbol\ell(\mathbf{p})$ dominates~it, noting that an evaluation point $\boldsymbol\ell(\mathbf{p}_\bullet)$ dominates $\boldsymbol\ell(\mathbf{p})$, denoted $\boldsymbol\ell(\mathbf{p}_\bullet) \preccur \boldsymbol\ell(\mathbf{p})$ or $\mathbf{p}_\bullet \preccur \mathbf{p}$, iff $\ell_i(\mathbf{p}_\bullet) \leq \ell_i(\mathbf{p})$ for all $i\in\{1,\ldots, n_\ell\}$ and $\ell_i(\mathbf{p}_\bullet) < \ell_i(\mathbf{p})$ for at least one~$i$.
\end{defn}
The optimization is contrasted to single-objective optimization where the solution correspond to one point (the~minimum).

\begin{remk}\label[remk]{th19}
	\Cref{th17} implies that a Pareto optimal point obtained using the simulation hardware stays Pareto optimal for the implementation hardware.
\end{remk}

In the proposed design problem, focus is given on the two underlying structural design parameters of the OCP, sampling time and prediction horizon. That is,
\begin{equation}
	\mathbf{p} \coloneqq \left(h,N\right) \in \mathcal {P}
	\coloneqq \mathbb{R}_{>0} \times \mathbb{N}_{>0}. \label{eq11}
\end{equation}
The rest of the coupled design parameters and all the decoupled design parameters are fixed. With such a focus, the resulting design problem would have specific characteristics, as revealed next, allowing for a specialized solver algorithm to be proposed.

\section{Key properties of the MOD-MPC system} \label{sc04}

This section analyzes the MOD-MPC system to reveal its key attributes. These include the smoothness properties and bounds on the objective and design parameters, establishing the underlying assumptions and guarantees that are useful for the subsequent development of a numerical solver. Based on the analytical foundation, a targeted solution method that is both effective and efficient can be appropriately developed.

\subsection{Smoothness of the design objectives}

\subsubsection*{Monotonicity}

The model used for the required computational resource~$\eta$~\eqref{eq09} has the following monotonicity property.

\begin{prop}[Monotonicity of $\eta$]
Consider the required computational resource $\eta(\cdot)$ in \eqref{eq09} and that \cref{th20} holds. For $\mathbf{p}\coloneqq(h,N)$, $\mathbf{p}\mapsto \eta(\mathbf{p})$ is:\\
$\bullet$ monotonically decreasing with respect to $h$,\\
$\bullet$ monotonically increasing with respect to $N$.
\label[prop]{th01}\end{prop}

\begin{proof}
	The monotonicity of $\eta$ can be directly taken from the dependence of $\eta$ to $h$ and $N$ as given in \eqref{eq09}.
\end{proof}

On the other hand, the control performance as measured by the closed-loop value function $V$ is non-monotonic.


\begin{theo}[Non-monotonicity of $V$ with respect to $h$]
Consider the control performance $V(\mathcal{X}_0,\cdot)$, $\mathcal{X}_0\subseteq\mathcal{X}_\text{\normalfont S}$ as in \eqref{eq10} and that \cref{th09,th24,th23,th22} hold. For $\mathbf{p}\coloneqq(h,N)$, $h\mapsto V\left(\mathcal{X}_0,\mathbf{p}\right)$ is monotonically increasing $\forall h\in\mathcal{H}_+$ and monotonically decreasing $\forall h\in\mathcal{H}_-$, therefore generally non-monotonic with respect to~$h$.
\label[theo]{th02}\end{theo}

\begin{proof}
Consider the closed-loop system\textsuperscript{\getrefnumber{fn01}}. For a fixed $N$, increasing $h$ both desirably increases the prediction horizon length $T\coloneqq Nh$ but also undesirably slowing down the sampling rate $1/h$ of the closed-loop system. At the limit $h\rightarrow0$, the system performs badly since the control prediction barely, if at all, captures any dynamics of the system. As $h$ is increased, performance is improved, until some point where all important dynamics are captured. Increasing prediction length is no longer as influential as the delayed sampling rate, after which the performance is worsened with increasing $h$. Therefore, the value function $U(\mathsf{x}_0,(h,N))$ is generally non-monotonic with respect to $h$. This is so that there are two mutually exclusive sets $\mathcal{H}_+$ and $\mathcal{H}_-$, where $\mathcal{H}_+\cup\mathcal{H}_-=\mathbb{R}_{>0}$. For all $h\in\mathcal{H}_+$, $\exists h_+> h$ such that the value function is increasing, $U(\mathsf{x}_0,(h_+,N))> U(\mathsf{x}_0,(h,N))$. $\mathcal{H}_-$ is defined similarly. $V$ is a linear combination of $U$ and thus has the same monotonicity properties.
\end{proof}

\stepcounter{footnote}
\footnotetext{with an OCP \eqref{eq01} satisfying \cref{th09,th24,th23,th22} with sampling period $h$, $N$ prediction steps, a solution $(x^*,u^*)$ and closed-loop value function $U(\mathsf{x}_0,(h,N))$, $\forall \mathsf{x}_0\in \mathcal{X}_\text{S}$ as in \eqref{eq02} \label{fn01}.}

\begin{theo}[Non-monotonicity of $V$ with respect to $N$]
Consider the control performance $V(\mathcal{X}_0,\cdot)$, $\mathcal{X}_0\subseteq\mathcal{X}_\text{\normalfont S}$ as in \eqref{eq10} and that \cref{th09,th24,th23,th22} hold. For $\mathbf{p}\coloneqq(h,N)$, $N\mapsto V\left(\mathcal{X}_0,\mathbf{p}\right)$ is monotonically increasing $\forall N\in\mathcal{N}_+$ and monotonically decreasing $\forall N\in\mathcal{N}_-$, therefore generally non-monotonic with respect to~$N$.
\label[theo]{th11}\end{theo}

\begin{proof}
Consider the closed-loop system\textsuperscript{\getrefnumber{fn01}}. For a fixed $h$, increasing $N$ increases the prediction horizon length $T\coloneqq Nh$. As a result more dynamics, as well as plant-model mismatch, are captured by the prediction. There is a trade-off balance between the two so that value function is generally non-monotonic with respect to $N$. This is so that there are two mutually exclusive sets $\mathcal{N}_+$ and $\mathcal{N}_-$, where $\mathcal{N}_+\cup\mathcal{N}_-=\mathbb{N}_{>0}$. For all $N\in\mathcal{N}_-$, $\exists N_-> N$ such that the value function is decreasing, $U(\mathsf{x}_0,(h,N_-))< U(\mathsf{x}_0,(h,N_+))$. $\mathcal{N}_+$ is defined similarly. $V$ is a linear combination of $U$ and thus has the same monotonicity properties. 
\end{proof}

\cref{th02,th11} are confirmed by the numerical observations in \cite{Bachtiar:15,Bachtiar:16}.


\subsubsection*{Continuity and differentiability}

The continuity of the solution of the OCP \eqref{eq01} with respect to $\mathbf{p}$ is described in the following.

\begin{lemm}[Continuity of $u^*$ \cite{Bachtiar:15,Bachtiar:16}] \label[lemm]{th10}
Consider the OCP~\eqref{eq01} satisfying \cref{th09,th24,th23,th22}. The unique optimal solution $\mathbf{p}\mapsto u^* (\mathsf{x}_i,\mathbf{p})$ of the OCP is differentiable with respect to $h$ for a given $N$.
\end{lemm}\vspace{-3mm}

\begin{proof}
The proof is given in Lemma~14 in \cite{Bachtiar:16}.
\end{proof}

\begin{theo}[Continuity of $V$]
Consider the control performance $V(\mathcal{X}_0,\mathbf{p})$, $\mathcal{X}_0\subseteq\mathcal{X}_\text{\normalfont S}$, and that \cref{th09,th24,th23,th22} hold. If \cref{th07} holds, $h\mapsto V\left(\mathbf{p}\right)$ is continuous with respect to $h$.
\label[theo]{th03}\end{theo}

\begin{proof}
	Let $\mathsf{f}(\cdot,\cdot)$ satisfy \cref{th07}. Let $\mathsf{z}(\cdot,\cdot)$ be the solution of $\dot{\mathsf{x}}=\mathsf{f}(\mathsf{x},\mathsf{u})$. $\mathsf{u}\mapsto\mathsf{z}(t,\mathsf{u})$ is continuous with respect to $\mathsf{u}$ (Theorem~3.5 in \cite{Khalil:02}). The control law \cref{eq02d} is given by $\mathsf{u}(\tau)=u_0^*(\mathsf{x}(t_i),\mathbf{p})$, $\forall \tau\in[ih,ih+h)$, $\forall i\in\mathbb{N}_{\geq0}$. Subsequently, from \cref{th10}, $\mathsf{u}$ is continuous and differentiable with respect to $h$, implying that the solution $h\mapsto\mathsf{z}(t,u^*(\mathsf{x}_i,(h,N)))$ is continuous with respect to $h$ for a given $N$. $h\mapsto U(\mathsf{x}_0,(h,N))$, $\forall \mathsf{x}_0\in \mathcal{X}_\text{S}$ is thus continuous with respect to $h$ for a given $N$ and satisfaction of \cref{th09,th24,th23,th22}. $V$ is a linear combination of $U$ and thus has the same monotonicity properties.
\end{proof}

The differentiability of the closed-loop value function can be described and is stated in the following.

\begin{theo}[Differentiability of $V$]
Consider the control performance $V(\mathcal{X}_0,\mathbf{p})$, $\mathcal{X}_0\subseteq\mathcal{X}_\text{\normalfont S}$, and that \cref{th09,th24,th23,th22} hold. If \cref{th07,th08} hold, $\mathbf{p}\mapsto V\left(\mathbf{p}\right)$ is differentiable with respect to $h$ for a given $N$.
\label[theo]{th04}\end{theo}

\begin{proof}
	Let $\mathsf{f}(\cdot,\cdot)$ satisfy \cref{th07,th08}. Let $\mathsf{z}(\cdot,\cdot)$ be the solution of $\dot{\mathsf{x}}=\mathsf{f}(\mathsf{x},\mathsf{u})$. $\mathsf{u}\mapsto\mathsf{z}(t,\mathsf{u})$ is differentiable with respect to $\mathsf{u}$ (Theorem~3.5 and Section~3.3 in \cite{Khalil:02}). The control law \cref{eq02d} is given by $\mathsf{u}(\tau)=u_0^*(\mathsf{x}(t_i),\mathbf{p})$, $\forall \tau\in[ih,ih+h)$, $\forall i\in\mathbb{N}_{\geq0}$. Subsequently, from \cref{th10}, $\mathsf{u}$ is continuous and differentiable with respect to $h$, implying that the solution $h\mapsto\mathsf{z}(t,u^*(\mathsf{x}_i,(h,N)))$ is differentiable with respect to $h$ for a given $N$. $h\mapsto U(\mathsf{x}_0,(h,N))$, $\forall \mathsf{x}_0\in \mathcal{X}_\text{S}$ is thus differentiable with respect to $h$ for a given $N$ and satisfaction of \cref{th09,th24,th23,th22}. $V$ is a linear combination of~$U$ and thus has the same monotonicity properties.
\end{proof}

\subsection{Competing nature of the design objectives}

The competing nature of a pair of functions that are both to be minimized (or maximized) is defined as follows.

\begin{defn}[Competing functions]
Two functions $a\mapsto f(a)$ and $a\mapsto g(a)$ are competing with each other for the design set $[\underline{a},\overline{a}]$ iff $a\mapsto f(a)$ is monotonically increasing and $a\mapsto g(a)$ is monotonically decreasing (or vice versa) on $[\underline{a},\overline{a}]$.
\label[defn]{th13}\end{defn}

In the multi-objective design of MPC, both the closed-loop value function $V$ and required computational resource $\eta$ are to be minimized. Based on \cref{th01,th02,th11}, the two design objectives of control performance and required computational resource are competing. This is detailed in the following.

\begin{lemm}[Competing design objectives]
The objective functions $V(\mathcal{X}_0,\cdot)$, $\mathcal{X}_0\subseteq\mathcal{X}_\text{\normalfont S}$, and $\eta(\cdot)$ are competing as per \cref{th13} within the design parameter set $\mathcal{P}_\text{\normalfont c} = (\mathcal{H}_+\times \mathcal{N}_-)$ from \cref{th01}, \cref{th02,th11}, and given that \cref{th09,th24,th23,th20,th22} hold.
\label[lemm]{th14}\end{lemm}

\begin{proof}
Consider \cref{th01}, \cref{th02,th11} and satisfaction of \cref{th09,th24,th23,th20,th22}. $\eta$ is monotonically decreasing with respect to $h$ and increasing with~$N$. $V$ is non-monotonic with respect to $h$ and $N$, and there exist a set such that $V$ is increasing with respect to $h$ and decreasing with respect to $N$, given by $\mathcal{P}_\text{\normalfont c}=(\mathcal{H}_+\times \mathcal{N}_-)$ from \cref{th02,th11}. Within this set,  $\forall\mathbf{p}\in\mathcal{P}_\text{\normalfont c}$, $\mathbf{p}\mapsto V(\mathbf{p})$ is monotonically increasing whilst $\mathbf{p}\mapsto \eta(\mathbf{p})$ is monotonically decreasing, or vice~versa.
\end{proof}

The search space $\mathcal{P}_\text{s}$ is assumed to intersect with $\mathcal{P}_\text{c}$ so that the MOD-MPC solution $\mathcal{P}_\bullet$ exist.
\begin{assm}
$\mathcal{P}_\bullet = \mathcal{P}_\text{\normalfont s}\cap\mathcal{P}_\text{\normalfont c}\neq\emptyset$.
\label[assm]{th21}\end{assm}
The associated Pareto front $\mathcal{L}$ as per \eqref{eq07} consists of Pareto optimal points as defined in \cref{th15}, each of which is a Pareto design choice  $\mathbf{p}\in\mathcal{P}_\bullet\subseteq\mathcal{P}_\text{s}$.

Further to its Pareto optimality, the quality of a point can be specified by its rank, as defined in the following.

\begin{defn}[Rank]
Given a countable set of points~$\mathcal{P}$, if a point $j$ is Pareto optimal then its rank $r_j=1$. Subsequently, a point $j$ has rank $r_j=\rho$ if it is Pareto optimal in $\mathcal{P}\setminus\mathcal{P}_{\rho-1}$ where $\mathcal{P}_\rho$ is the set of all points with rank $r\leq\rho$.
\label[defn]{th06}\end{defn}

That is, all Pareto optimal points in a given set of points have a rank of 1. The Pareto optimal points in the set that excludes points with rank 1 have rank 2, and so on.




\subsection{Bounds on the Pareto design set}

The first bounding of the design parameter comes from the fact that it is  numerically impractical to search the open set $\mathcal{P}$ in~\eqref{eq11}. Hence, the search space $\mathcal{P}_\text{s}$ in~\eqref{eq03} must be a closed set $\mathcal{P}_\text{s} \subset \mathcal{P}$ that is able to be practically searched to find the Pareto design set $\mathcal{P}_\bullet$.

Next, an assumption on the Pareto design set can be made based on some intuitions on the nature of the design problem.
\begin{assm}[upper-bound on $h$]
	For a given $N>1$, the Pareto design set is upper-bounded by $\hat{h}$. This bound is defined by the notion that $\exists \hat{h}\in \mathcal{P}_\text{\normalfont s}$ such that $\forall h>\hat{h}$, $\boldsymbol\ell((h,N-1)) \preccur \boldsymbol\ell((h,N))$.
\label[assm]{th16}\end{assm}
The assumption comes from the fact that as the sampling period $h$ is increased for a given number of prediction steps $N$, the competitive effect of reducing the required computing resource will be diminished and surpassed by the deterioration in control performance as a result of more infrequent sampling. At this point, using fewer prediction steps, e.g.\ $N-1$, would reduce resource requirement more competitively than deteriorating performance. The opposite is also true; as the sampling period is shortened, the effect of increasing control performance due to a finer sampling would eventually be overtaken by an increase in computing resource required, after which increasing~$N$, to e.g.\ $N+1$, would be more competitive in increasing performance.

\begin{assm}[Lower bound on $h$]
	For a given $N\geq1$, the Pareto design set is lower bounded by $\check{h}$. This bound is defined by the notion that  $\exists \check{h}\in \mathcal{P}_\text{\normalfont s}$ such that $\forall h<\check{h}$, $\boldsymbol\ell((h,N+1)) \preccur \boldsymbol\ell((h,N))$.
\label[assm]{th12}\end{assm}

Based on the bounds on the sampling period and the smoothness properties of the design objectives defined in the earlier subsection, bounds can be specified for $\mathcal{P}_\bullet$.

\begin{theo}[Bound on $\mathcal{P}_\bullet$]\label[theo]{th05}
Consider a rectangular search space $\mathcal{P}_\text{\normalfont s}=\big\{(h,N):h\in[\underline{h},\overline{h}]\text{ and }N\in[\underline{N},\overline{N}]\big\}$ and that \cref{th16,th12} is satisfied. The Pareto design set $\mathcal{P}_\bullet$ is contained within the band $\mathcal{P}_\text{\normalfont b} \coloneqq \left\{\mathbf{p} : h\leq m_1N + h_1 \text{ and } h\geq m_2N + h_2\right\}$ for some negative gradients $m_i\in\mathbb{R}_{<0}$ and constants $h_i$, $i\in\{1,2\}$.
\end{theo}

\begin{proof}
Consider a Pareto design set $\mathcal{P}_\bullet$ in a rectangular search space $\mathcal{P}_\text{\normalfont s}=\big\{(h,N):h\in[\underline{h},\overline{h}]\text{ and }N\in[\underline{N},\overline{N}]\big\}$. For a given number of prediction steps $N>1$, let the sampling periods corresponding to the Pareto design set have an upper-bound from \cref{th16} denoted $\hat{h}_N$. Paraphrasing the assumption, $\exists{h}$ such that $\boldsymbol\ell((h,N-1)) \preccur \boldsymbol\ell((h,N))$, $\forall h>\hat{h}_N$ for the given upper-bound $\hat{h}_N$. Consequently, there must exist an upper-bound associated with $N-1$ prediction steps $\hat{h}_{N-1}$ that is larger than the upper-bound $\hat{h}_{N}$, giving $\hat{h}_N < \hat{h}_{N-1}$. Therefore, the Pareto design set can be upper-bounded by a line of negative gradient with respect to $N$. An opposite notion can be made using \cref{th12} to form a lower bound with a negative gradient. This gives a bound in the form of a band $\mathcal{P}_\text{\normalfont b}$ as in the theorem.
\end{proof}

\section{Numerical solution to the MOD-MPC problem} \label{sc05}

\subsection{Effective and efficient solver characteristics}

Analysis of the key properties of the system results in a number of characteristics required by a proposed numerical optimizer used to solve the MOD-MPC problem~\eqref{eq03} accurately and quickly, as summarized below.

\begin{cond}[Convergent] \label[cond]{cn01}
	\Cref{th14} implies that there is a Pareto optimal design set $\mathcal{P}_\bullet$ for a given search space $\mathcal{P}_\text{s}$ associated with the trade-off of the competing objectives. The solver should be able to effectively find $\mathcal{P}_\bullet$ with certain guarantees.
\end{cond}

\begin{ccond}{1\text{a}}[Global] \label[cond]{cn01a}
	\Cref{th01,th02,th11} define that the required computational resource $\eta$ is monotonic and that the value function $V$ is non-monotonic. The solver needs to search globally and handle the many local optima on the objective surface.
\end{ccond}

\begin{ccond}{1\text{b}}[Able to handle discrete parameters] \label[cond]{cn01b}
	The solver must be able to handle discrete design parameters defined in \eqref{eq11}.
\end{ccond}

\begin{remk}[]\label[remk]{th25}
	\Cref{cn01a,cn01b} are necessary conditions for \cref{cn01} to be fulfilled. Satisfaction of these two is not always sufficient to satisfy \cref{cn01}.
\end{remk}

\Cref{cn01} (necessarily with \cref{cn01a,cn01b}) is a sufficient condition for a numerical solver to be accurate (convergent) for the MOD-MPC problem. Additional features are necessary for the solver to converge quickly and efficiently.

\begin{cond}[Continuous] \label[cond]{cn02}
	The solver could rely on \cref{th03,th04} that define the continuity/differentiability of the value function $V$ based on the knowledge of $\mathsf{f}(\cdot,\cdot)$.
\end{cond}

\begin{cond}[Focused] \label[cond]{cn03}
	\Cref{th05} states that for a rectangular search space, the Pareto optimal design set is located within a specific space defined as a band. As a consequence, the solver should be able to focus its search within the band and omit any ineffectual space.
\end{cond}

\Cref{cn02,cn03} are sufficient for the solver to be efficient and performs better than a general-purpose solver.


\subsection{A compliant solver algorithm (DITRI)}

Having now specified the sufficient conditions for an effective and efficient solver, a specialized solver satisfying all the conditions can be proposed for the MOD-MPC problem. An algorithm is proposed based on Lipschitzian optimization \cite{Jones:93}, denoted `DIviding TRIangles' (DITRI), with details outlined in \ref{apA}.

The choice of a Lipschitzian approach addresses \cref{cn01a,cn01b} necessary for the solver to be convergent. Lipschitzian optimization is gradient-free and is built for a global search, therefore addressing \cref{cn01a}. The method directly handles discrete design parameters specified in \cref{cn01b}. As a whole, the proposed solver is guaranteed to converge and satisfies \cref{cn01}, making it an accurate solver that can effectively find the solution (Pareto design set) $\mathcal{P}_\bullet$ of the MOD-MPC problem \eqref{eq03} in a given search space.

\begin{theo}[Convergent DITRI]\label[theo]{th18}
	Consider a multi-objective optimal design problem with objectives $\mathbf{p}\mapsto\boldsymbol\ell(\mathbf{p})$ for $\mathbf{p}\in\mathcal{P}_\text{\normalfont s}$ where $\mathcal{P}_\text{\normalfont s}$ is a finite search space. Let the solution of the problem be $\mathcal{P}_\bullet\subseteq\mathcal{P}_\text{\normalfont s}$. Also let $\{\mathcal{S}_i\}_{i=0}^{\overline{i}}$ for some indexing variable $i$ be a sequence of solutions generated from a given initialization $\mathcal{S}_0$. DITRI, described in Algorithm~\ref{al01}, is a convergent algorithm such that $\lim_{i\rightarrow i^\text{\normalfont c}} S_i=\mathcal{P}_\bullet$ for some $i^\text{\normalfont c}<\infty$.
\end{theo}

\begin{proof}
	The proof is given in \ref{apC}.
\end{proof}

As well as being accurate, DITRI is designed to be an efficient optimizer by using projection of bounds (\ref{apA-01}). This assumes continuity in the objective function and satisfies \cref{cn02}. Finally, DITRI conducts a focused search of a given search space as outlined in \ref{apA-02}--4 to fulfill \cref{cn03}.

\section{Simulation results} \label{sc06}

Two real-world examples are investigated. $Q$ and $R$ are chosen accordingly for each case. $Q_\text{f}=P_\text{ARE}$, where $P_\text{ARE}$ is the solution of the algebraic Riccati equation for the simulated plant. The OCP is represented as a sparse QP. The global search criterion \eqref{eq12} is set as $\overline{d}(i)\coloneqq \sqrt{5/9^{i/8}}$ so that no search space is larger than that equivalent to $i/8$ divisions from the initial triangle.

\begin{figure}[t!]
\centering
\includegraphics[scale=0.75,trim=0 15 0 0]{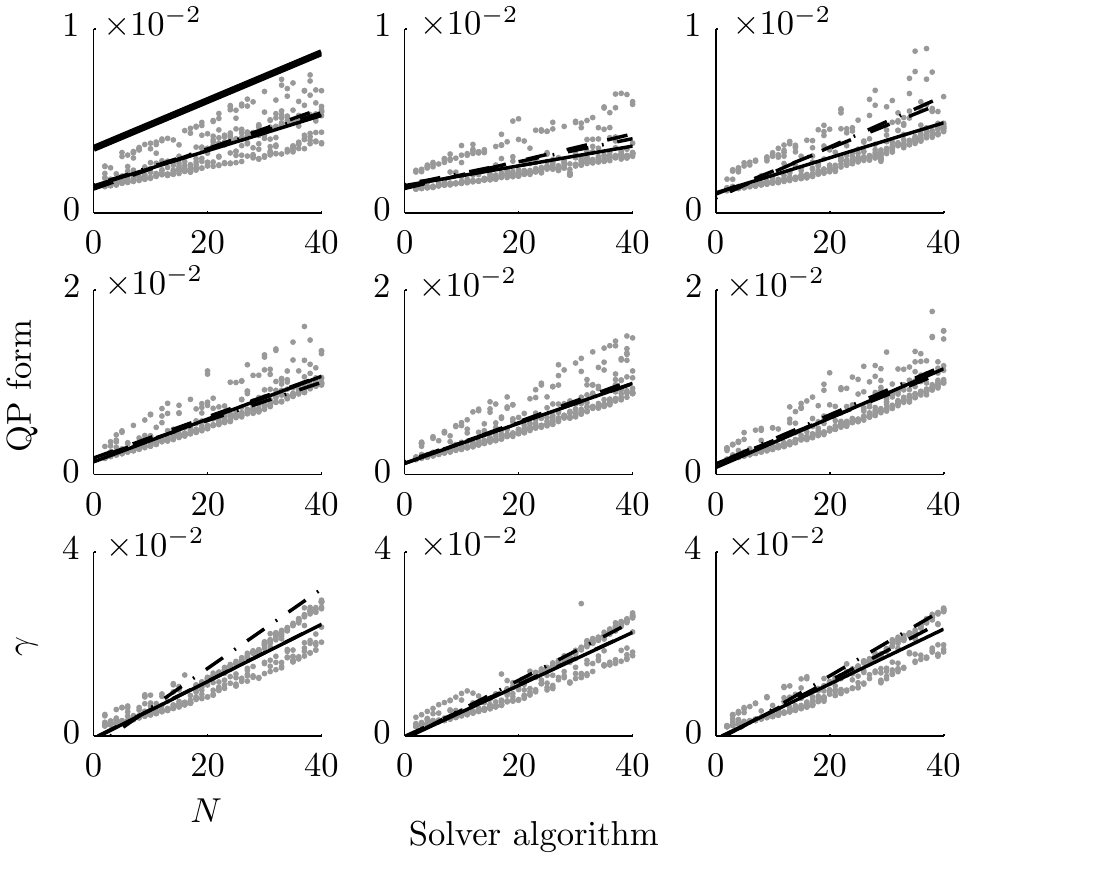}
\caption[]{Relation between solution time and $N$ for different QP representations and QP solver algorithms. The upper-bound model $\gamma$ used in this paper is the thick solid line on the top-left graph. Left-to-right: Interior point, primal simplex and concurrent methods. Top-to-bottom: sparse, sparse-delta \cite{Longo:14} and dense QP representations. Data is obtained from the PAA problem \eqref{eq13} given in Subsection~\ref{sc06.1}. Across all plots, gray point plots are data for $h=5$ ms (point plots for other $h$ values are not shown). The solid, dashed and dot-dashed thin black lines are linear fits for the data with $h=5$ ms, $10$ ms and $5$ s.}
\label{fig02}
\end{figure}

The simulation hardware specifications are given in Table~\ref{tb02}. Gurobi \cite{Gurobi:XX} in MATLAB is used in this study as the OCP solver. The algorithm used is the interior point (barrier) method, with all tolerances set as the default. A representative result for the relationship between the solution time and number of prediction steps $N$ for a range of sampling time $h$ is shown in Fig.~\ref{fig02}. It is shown that $\gamma$ is generally increasing with $N$ and that it is very weakly correlated to $h$, verifying \cref{th20}. From the obtained data, the relationship for the chosen algorithm and QP form is mostly linear. $\gamma$ in \eqref{eq16} is modeled as a linear function with the chosen constants $a_1=1.3\times10^{-4}$ and $a_0=3.5\times10^{-3}$. The relevant data and model is shown on the top-left graph of the figure.

\begin{table}[t!]
\setlength\tabcolsep{1.5mm}
\centering\small
\begin{tabular}{rl}
\toprule
{Cores}                & 4           \\
{Cycle frequency/core} & 3.4-3.9 GHz \\
{Operations/cycle}     & 8           \\
{FLOP/s}               & $\sim$$109$-$125\times10^9$\\
\midrule
{Cores available for simulation}                & 1           \\
{FLOP/s available for simulation}               & $\sim$$27$-$31\times10^9$\\ \bottomrule
\end{tabular}
\caption{Simulation hardware specifications (Intel\textsuperscript{\textregistered} Core\texttrademark i7-3770 Processor) \cite{Intel:XX}.}
\label{tb02}
\end{table}

\subsection[]
{Test plant models\footnote{Symbols used in a plant model are used exclusively in the model and should not be confused with symbols introduced elsewhere.}} \label{sc06.1}


The first test case looks at a missile pitch-axis autopilot (PAA) adapted from \cite{Bachtiar:14}. The missile is flown at a cruising altitude and the autopilot is to control the missile to track a commanded acceleration. The second case aims to design a controller for diesel engine control. The engine is modeled by a mean-value engine model (MVEM) taken from \cite{Broomhead:15a} and the control objective is to track a given engine speed and power output.

A coordinate shift is applied appropriately to transform a given tracking problem into that of regulation,
\begin{equation}
\begin{aligned}
  &\mathsf{x} = \mathsf{e} - \mathsf{e}_\text{s} \hspace{6mm}
  \mathsf{x}_0 = \mathsf{e}_0 - \mathsf{e}_\text{s} \hspace{6mm}
  \overline{x} = \overline{e} - \mathsf{e}_\text{s} \\
  &\mathsf{f}_\text{e}(\mathsf{x}+\mathsf{e}_\text{s},\mathsf{u}) = \mathsf{f}(\mathsf{x},\mathsf{u})
  \hspace{15mm} \text{etc.}
\end{aligned}\label{eq18}
\end{equation}

\subsubsection*{Pitch-axis autopilot}

The first test case looks at a missile pitch-axis autopilot (PAA) at $20\,000$ ft. The nonlinear tracking model is
\begin{equation}\label{eq13}
	\mathsf{f}_\text{e}(\mathsf{e},\mathsf{u}) =
		\begin{bmatrix*}[l]
			\mathsf{e}_2 + \cos(\mathsf{e}_1)F(\mathsf{e}_1,\mathsf{u})/(mv) \\
			L(\mathsf{e}_1,\mathsf{u})/I_y \\
			\mathsf{e}_4 \\
			-\omega_0^2 \mathsf{e}_3 - 2\zeta\omega_0 \mathsf{e}_4 + \omega_0^2 \mathsf{e}_5\\
			\mathsf{u}
		\end{bmatrix*}
\end{equation}
where $F$ and $L$ are the nonlinear mapping for the aerodynamic lift force and pitching moment respectively. $\mathsf{e}_1$ is the angle of attack and $\mathsf{e}_2$ is the pitch rate of the missile. The actuation of the fin deflection $\mathsf{e}_3$ is modeled as a second order system. The input $\mathsf{u}$ is the rate of the commanded fin deflection $\mathsf{e}_5$. Missile speed $v=Mv_s$, where $v_s$ is the speed of sound at $20\,000$ ft, is constant at Mach number $M=2.5$. $m$ and $I_y$ are the mass and moment of inertia of the missile respectively. These parameter values, other missile frame parameters, constants related to the actuation dynamics, along with the aerodynamic coefficients and models used for $F$ and $L$ are the same as given in \cite{Bachtiar:14,Nichols:93}.



The control objective is to track a given acceleration output $\mathsf{y}=F/(mg)$, where $g$ is the gravitational acceleration. The test scenario is to track 5 different acceleration outputs from steady-state at $0g$, $\mathsf{e}_0=0$. The outputs are $2$, $4$, $6$, $8$ and $10g$, each associated with a unique steady state $\mathsf{e}_\text{s}$ and initial condition as per~\eqref{eq18}, making up the set of initial conditions $\mathcal{X}_0$ in~\eqref{eq10} that are equally weighed, $w_i=0.2$, $i\in\{0,\ldots,5\}$. The states are upper- and lower- bounded by $\overline{e}=(20^\circ, \; 35^\circ/\text{s}, \; 45^\circ, \; 10^{6\,\circ}/\text{s}, \; 45^\circ)$ and $\underline{e}=-\overline{e}$. The input is bounded by $\overline{u} = 10^{6\,\circ}/\text{s}$ and $\underline{u}=-\overline{u}$. $Q=C^\mathsf{T}C$ where $C\coloneqq  d\mathsf{y}/d\mathsf{e} \big|_{\mathsf{e}_\text{s},0}$ comes from the linearization of the output $\mathsf{y}$ at the target steady state and $R=10^{-6}$. 

\subsubsection*{Diesel engine control with a mean-value engine model}

The second case looks at engine control with a 5-state 3-input mean-value engine model (MVEM) from \cite{Broomhead:15a}. The 5 states are the engine speed $\mathsf{e}_1$, turbine speed $\mathsf{e}_2$, VGT actuator position $\mathsf{e}_3$, intake manifold pressure $\mathsf{e}_4$ and temperature $\mathsf{e}_5$. The three inputs are the injection duration $\mathsf{u}_1$, load applied to the engine by the generator $\mathsf{u}_2$ and the VGT commanded position $\mathsf{u}_3$. The model is
\begin{equation}\label{eq08}
	\mathsf{f}_\text{e}(\mathsf{e},\mathsf{u}) =
		\begin{bmatrix*}[l]
			(\tau_\text{eng}-\mathsf{u}_2)/J_\text{e} \\
			(P_\text{t}-P_\text{c})/(J_\text{t}\mathsf{e}_2) \\
			(\mathsf{u}_3-\mathsf{e}_3)/\tau_\text{VGT} \\
			\dfrac{R_\text{a}}{V_\text{im}}(\dot{m}_\text{c} + \dot{m}_\text{EGR} - \dot{m}_\text{ei}) \mathsf{e}_5 
			 + \dfrac{\mathsf{e}_4}{\mathsf{e}_5}\dot{\mathsf{e}}_5 \\[2.5mm]
			\dfrac{R_\text{a}}{V_\text{im}c_\text{va}} \dfrac{\mathsf{e}_5}{\mathsf{e}_4} \Big(R_\text{a}(T_\text{ic}\dot{m}_\text{c} + T_\text{EGR} \dot{m}_\text{EGR} - \mathsf{e}_5 \dot{m}_\text{ei}) + \\
			\hspace{5mm} c_\text{va}\dot{m}_\text{c}(T_\text{ic}-\mathsf{e}_5) + c_\text{va}\dot{m}_\text{EGR}(T_\text{EGR}-\mathsf{e}_5) \Big)
		\end{bmatrix*}
\end{equation}
with static states determined after time-scale separation
\begin{equation}
	0 =
	\begin{bmatrix*}[l]
		\dot{m}_\text{cyl}(\mathsf{u}_1) - \dot{m}_\text{EGR} - \dot{m}_\text{t} \\
		O_\text{cyl}(\mathsf{u}_1)-O_\text{em} \\
		T_\text{em}-\dfrac{c_\text{pe}\dot{m}_\text{cyl}T_\text{cyl}(\mathsf{u}_1) + hA_\text{em,i}T_\text{em,s}}{c_\text{pe}\dot{m}_\text{cyl}(\mathsf{u}_1) + hA_\text{em,i}}
	\end{bmatrix*}
\end{equation}
and the further assumptions that
\begin{equation}
	0 = \begin{bmatrix*}[l] \mathsf{e}_\text{EGR} \\ O_\text{im} - O_\text{FRs}\end{bmatrix*}.
\end{equation}
$J_\text{e}$, $J_\text{t}$, $\tau_\text{VGT}$ and $V_\text{im}$ represent physical engine parameters. $R_\text{a}$ is the specific gas constant for the ambient gas. $c_\text{va}$ is the isometric specific heat of the ambience and $c_\text{pe}$ is the isobaric specific heat of the exhaust gas. $O_\text{FRs}$ is the stoichiometric mass ratio.

$\tau_\text{eng}$ is the engine load. $P$, $\dot{m}$ and $T$ denote power outputs, mass flows and temperatures respectively. Subscripts \{t\}, \{c\}, \{EGR\} and \{ic\} represent associations with the turbine, compressor, EGR and intercooler/compressor respectively. Subscripts \{em\}, \{em,s\} and \{em,i\} represent associations with the exhaust manifold. Subscripts \{cyl\} and \{ei\} represent associations with the cylinders of the engine. Expressions for these algebraic variables are given in \cite{Broomhead:15a}.

The initial state of the engine is at $2\,000$ rpm producing $20$ kW of power. The control objective is for the engine to track a steady-state at $2\,500$ rpm producing $36$ kW of power. The states are constrained with an upper-bound of $\overline{e}=\big(2500\,\text{rpm}, \; 150\,000\,\text{rpm}, \; 10^8, \; 10^8, \; 10^8 \big)$ and lower bound of $\underline{e}=(1500\,\text{rpm}, \; 45\,000\,\text{rpm}, \; 0, \; 0, \; 0 )$. The input is bounded by $\overline{u} = (1\,\text{ms}, \; 300\,\text{Nm}, \; 90)$ and $\underline{u}=( 0.5\,\text{ms}, \; 10\,\text{Nm}, \; 60 )$ respectively. $Q = \text{diag}(1,0,0,0,0)$ and $R=\text{diag}(0,1,1)$ as to track engine speed, power output, and VGT position. 

\subsection{Test results}

\begin{figure}[b!]
\centering
\includegraphics[scale=0.75,trim=0 0 0 10]{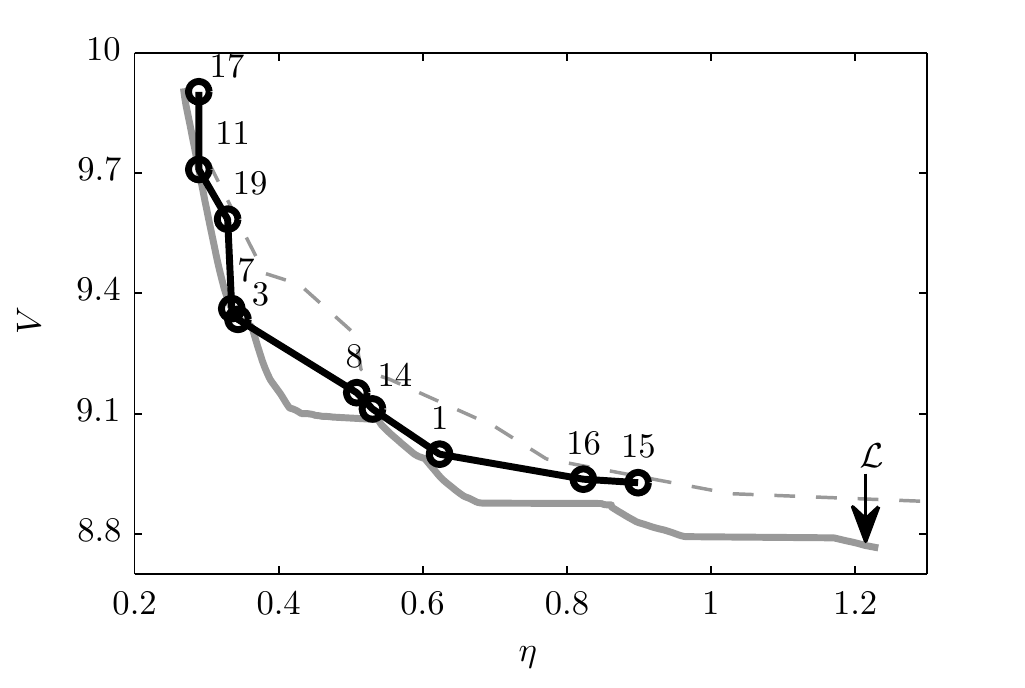}
\includegraphics[scale=0.75,trim=0 10 0 0]{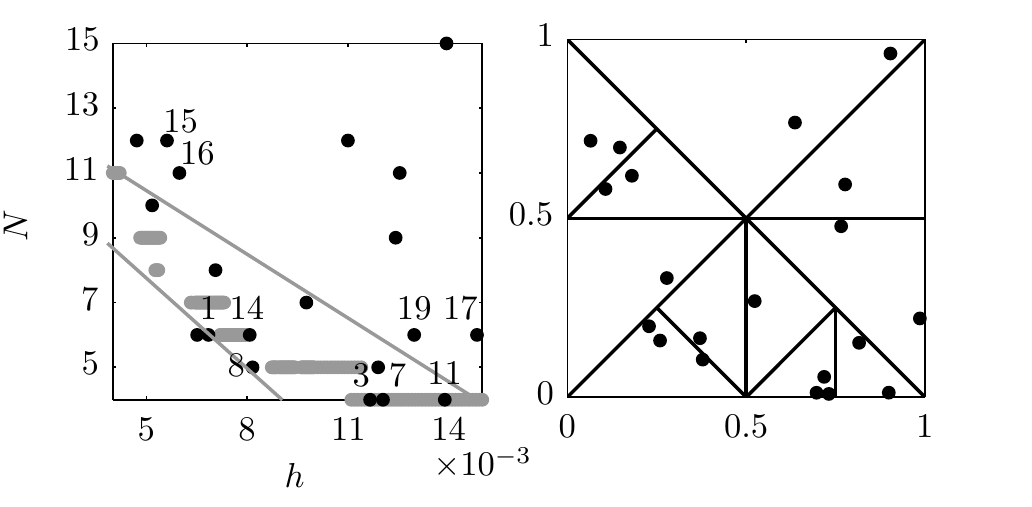}
\caption[]{Illustrative result for the PAA test case\textsuperscript{\getrefnumber{fn03}}.}
\label{fig03}
\end{figure}

\begin{figure}[b!]
\centering
\includegraphics[scale=0.75,trim=0 0 0 10]{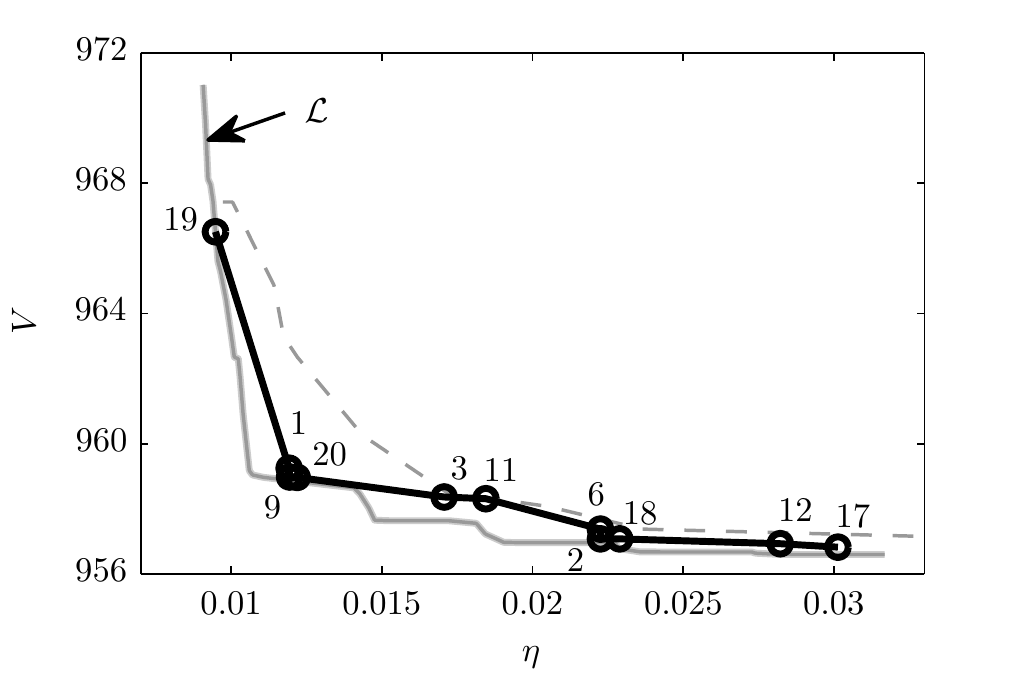}
\includegraphics[scale=0.75,trim=0 10 0 0]{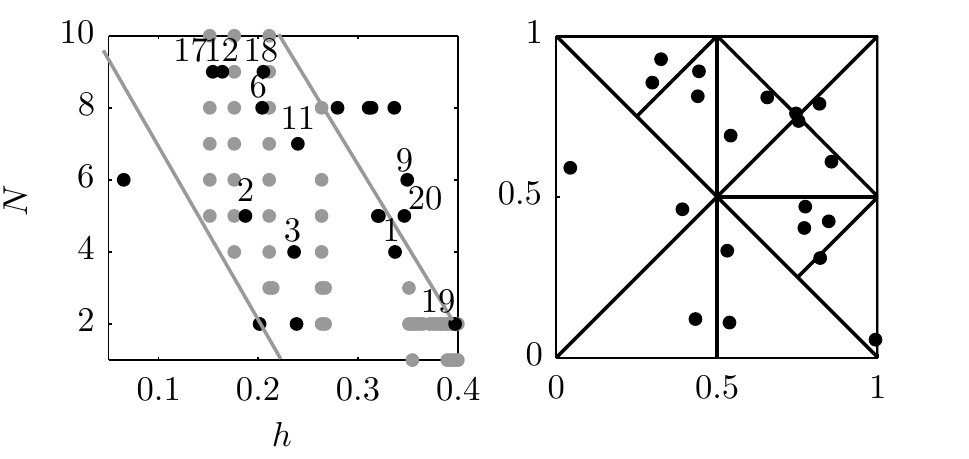}
\caption[]{Illustrative result for the MVEM test case\textsuperscript{\getrefnumber{fn03}}}
\label{fig04}
\end{figure}

Fig.~\ref{fig03} shows a representative result for the PAA case~\eqref{eq13}. A resulting trade-off curve is obtained after 20 evaluations using DITRI with $\overline{\mathbf{p}}=(0.015,\,15)$ and $\underline{\mathbf{p}}=(0.001,\,3)$, consisting of 10 different designs.  For the MVEM case~\eqref{eq08}, a representative result is shown on Fig.~\ref{fig04} for 20 evaluations using DITRI with $\overline{\mathbf{p}}=(0.4,\,10)$ and $\underline{\mathbf{p}}=(0.05,\,1)$. After 20 evaluations, 11 designs on a trade-off curve are obtained. The associated solution in the design parameter space is shown on the bottom plots of each figure.

The trade-off curves represent the set of optimal designs a practitioner can choose from. For example, design 16 of the PAA case (Fig.~\ref{fig03}) has a controller design with a sampling time of $h=6.5$ ms and $N=11$ prediction steps (prediction horizon of $\sim$70 ms). Consequently, the implementation hardware of the controller should be able to solve the OCP with $N=11$ in under 6.5 ms.

To help design the implementation hardware, the Resource Number of the design choice can be examined. PAA design 16 is associated with $\eta=0.82$, indicating that the implementation hardware must have at least 0.82 times the processing power in FLOP/s of the simulation hardware (Table~\ref{tb02}). This depends on hardware capabilities and implementation architecture, including clock-frequency, pipelining and parallel-processing.

The trade-off curves also reveal the sensitivity of control performance to computational resource. In the PAA case, performance improvement after $\eta\simeq0.6$ is not significant anymore. This implies that there is not much benefit to be gained from hardware more powerful than $\eta\simeq0.6$. In the MVEM case, the value is $\eta\simeq0.015$.

\stepcounter{footnote}
\footnotetext{Top: trade-off curve (black) along with the Pareto front $\mathcal{L}$ obtained from a full design exploration (gray). Each point is labeled by the associated evaluation number. The dashed-gray line is a non-convergent HVOL solution. Bottom-left: the associated solution in parameter space along with the band described in \cref{th05}. Bottom-right: accompanying plot showing the normalized space $c_1$-$c_2$ and triangle divisions used internally in DITRI. Total number of evaluations is $\overline{i^\text{ev}} = 20$.\label{fn03}}

\subsection{Validation of the prescribed solver characteristics}

To show the importance of the conditions for an effective and efficient solver prescribed in Section~\ref{sc05}, DITRI is compared to two other algorithms. The first is a non-dominated sorting genetic algorithm (NSGA) adapted from \cite{Deb:02} and the second is an algorithm based on surrogate hyper-volume improvement (HVOL) adapted from \cite{Tesch:13}. Table~\ref{tb03} outlines how each solver satisfies the specified conditions.

\begin{table}[b!]
\setlength\tabcolsep{1.5mm}
\centering\small
\begin{tabular}{rccccc}
\toprule
          & \multicolumn{5}{c}{Condition}                                  \\
Algorithm & 1          & 1a         & 1b         & 2          & 3          \\
\midrule
DITRI     & \checkmark & \checkmark & \checkmark & \checkmark & \checkmark \\
NSGA      & \checkmark & \checkmark & \checkmark &            &            \\
HVOL      &            & \checkmark & \checkmark &            &            \\ \bottomrule
\end{tabular}
\caption{Fulfillment of conditions for convergence and efficiency.}
\label{tb03}
\end{table}

\subsubsection*{Effective convergence (\cref{cn01})}

To assess the convergence of the trade-off curves obtained, the curves are compared to the true Pareto front~$\mathcal{L}$. Since the true Pareto front is not known, it is approximated by doing a full exploration on a uniform grid of $400$ $h$-values for each $N$-value in the parameter space. For both case studies, it is shown that the solution obtained by DITRI is close to the true Pareto front of the problem (Figs.~\ref{fig03} and \ref{fig04}). 

\begin{figure}[t!]
\centering
\includegraphics[scale=0.75,trim=0 15 0 5]{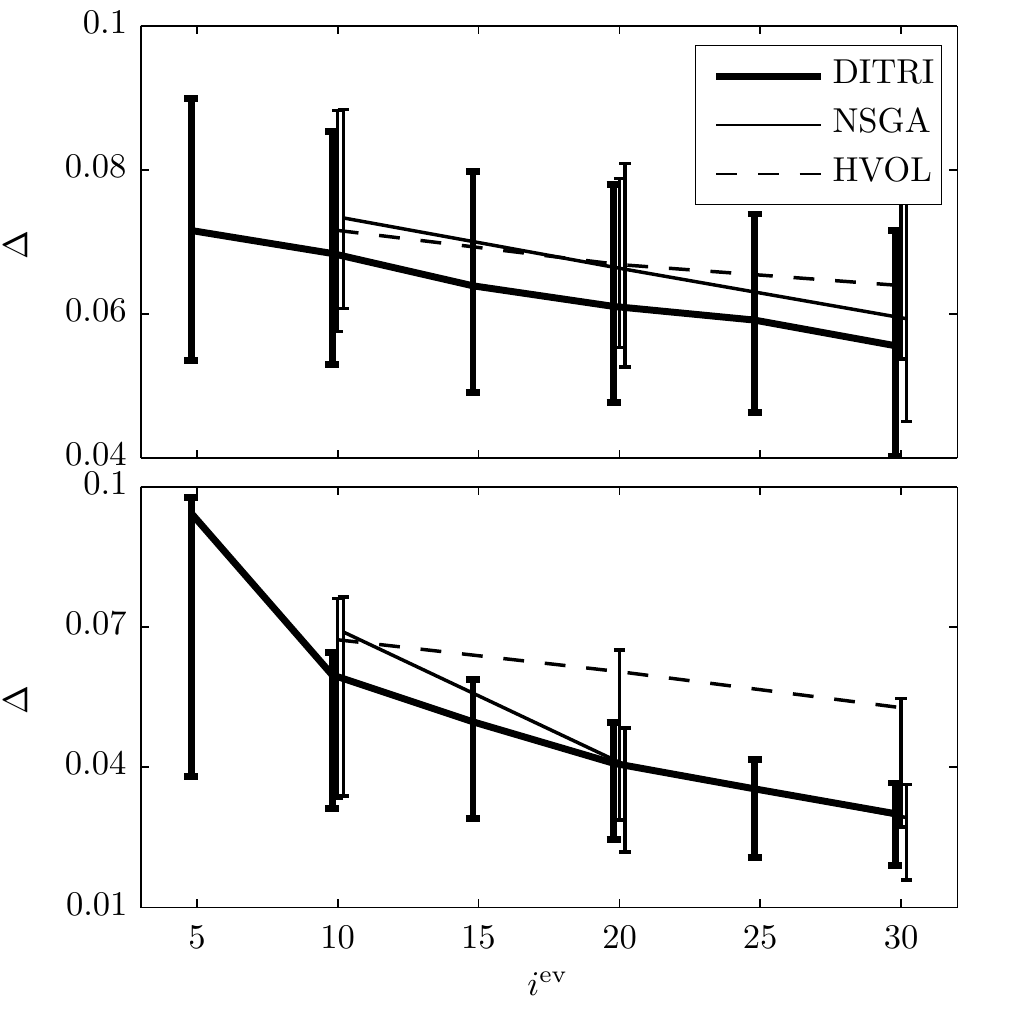}
\caption[]{Plot of closeness against evaluation count for the PAA (top) and MVEM (bottom) test cases\textsuperscript{\getrefnumber{fn04}}.}
\label{fig07}
\end{figure}

\begin{figure}[t!]
\centering
\includegraphics[scale=0.75,trim=0 15 0 5]{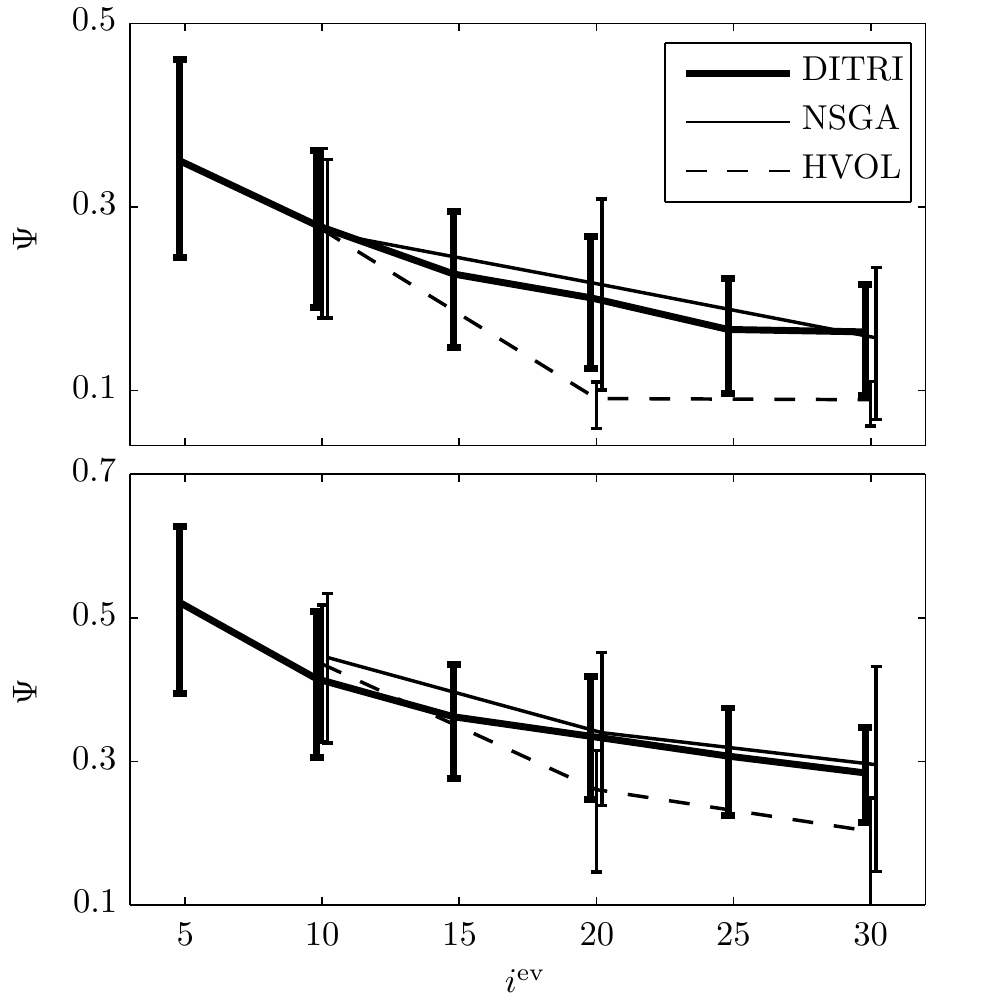}
\caption[]{Plot of coverage against evaluation count for the PAA (top) and MVEM (bottom) test cases\textsuperscript{\getrefnumber{fn04}}.}
\label{fig08}
\end{figure}

The closeness of a trade-off curve to the true Pareto front can be measured by calculating the average of the closest Euclidean distance between each point on the trade-off curve to the Pareto front. This measure is denoted $\Delta$ and plotted in Fig.~\ref{fig07} against function evaluation count. A second metric calculates the Euclidean distance of the tips (vertices) of the trade-off curve and Pareto front. This measures the coverage of the solution, denoted $\Psi$, and is shown in Fig.~\ref{fig08}. Calculation of both metrics are based on a normalized design objective values.

Figs.~\ref{fig07} and \ref{fig08} show that trade-off curves produced by DITRI and NSGA approach the Pareto front with increasing function evaluation counts. However, HVOL struggles to converge. Figs.~\ref{fig03} and \ref{fig04} show non-convergent trade-off curves, each from 100 evaluations using HVOL, confirming HVOL's inability to find the Pareto front. This is consistent with the expectation, since both DITRI and NSGA satisfies \cref{cn01} for convergence, whereas HVOL does not. The fulfillment of \cref{cn01a,cn01b} by HVOL is not sufficient to guarantee convergence, consistent with \cref{th25}.

\stepcounter{footnote}
\footnotetext{ The graphs show the mean across 500 trials with error bars showing the 25\textsuperscript{th} and 75\textsuperscript{th} percentiles.\label{fn04}}

\subsubsection*{Efficient search (\cref{cn02,cn03})}

The results in Figs.~\ref{fig03} and \ref{fig04} show that the Pareto solution lies within a band as described in \cref{th05}. DITRI takes advantage of this, unlike NSGA and HVOL. Comparing the results in Figs.~\ref{fig07} and \ref{fig08}, DITRI exhibits the best convergence rate. This is consistent with the fact that DITRI satisfies \cref{cn02,cn03} for efficiency, while the general-purpose NSGA and HVOL do not.

\section{Conclusions and future work} \label{sc07}\vspace{-3mm}

This paper presented an MPC design approach in a multi-objective fashion, treating control performance and the required computational resource as concurrent objectives in a given control problem. Focus was given to tuning the structural attributes of the MPC, namely the sampling time and prediction horizon. This approach is more comprehensive than those that explore only one design objective and treat software and hardware separately. A co-design of both MPC algorithm and hardware streamlines the design process, avoiding unnecessary costs. The proposed approach was studied analytically to present several theoretical results that reveal key properties of the design problem and subsequently prescribe necessary and sufficient conditions for an effective solver. Finally, two tests on real-world examples were conducted to demonstrate the design approach, as well as the importance of the conditions specified for an effective solver of the design problem.

Future work following the study could consider other coupled design parameters currently kept constant, such as the choice of prediction model. Furthermore, the scope can be extended further beyond the MPC structure to include attributes of the numerical method  used to solve the OCP, such as the algorithm and its tolerances, as well as features of the implementation hardware resource such as data representation type. This would extend the idea to a full co-design approach that looks at attributes of both software and hardware instead of focusing only on software parameters.

\appendix
\section{DITRI Algorithm} \label{apA}

\subsection{Projection of bounds} \label{apA-01}

The principles of Lipschitzian optimization are outlined in Algorithm~\ref{al02}. In each iteration, given a set of point(s)~$P$, a point $\mathbf{p}_j \in P$ is \textit{potentially optimal} if its \textit{projected} lower bound of the (minimized) objective within the associated search space must be equal or better than all points in $P$. The search space of all potentially optimal points will be partitioned into smaller divisions and a new point is the evaluated in each division.

\begin{algorithm}[b]
\caption{Generic Lipschitzian optimization}
\begin{algorithmic}[1]
\Require Search space bounds
\State Evaluate initial point(s)
\Repeat
	\State Find a set $O$ of potentially optimal points
	\ForAll {$o \in O$}
		\State \parbox[t]{\dimexpr\linewidth-2em}{Evaluate new points based on the search space division of potentially optimal point $o$}
	\EndFor
\Until iteration or evaluation count limit is reached
\end{algorithmic}
\label{al02}
\end{algorithm}

The bound \textit{projection} must be consistent throughout, dictated by a constant referred to as the Lipschitz constant, hence the name of the algorithm. The projected bound $L(d_j)$ of an evaluated objective $\ell(\mathbf{p}_j)$ in a search space of size $d_j$ is defined as
$
	L(d_j) = \ell(\mathbf{p}_j) - K_\text{L}d_j
$
for some Lipschitz constant $K_\text{L}>0$. For the projection to make sense, it assumes that the objective is continuous, as given in \cref{cn02}. For a point $\mathbf{p}_j$, the bound is better (smaller) if it has a smaller objective $\ell(\mathbf{p}_j)$ and/or a bigger search space size $d_j$. The point $j$ is potentially optimal if $L(d_j) \leq L(d_i)$ for all $i\in P$, that is
\begin{equation}
	\ell(\mathbf{p}_j) - K_\text{L}d_j \leq \ell(\mathbf{p}_i) - K_\text{L}d_i, \; \forall i\in P \text{ for a given } K_\text{L}>0.
\label{eq05}\end{equation}

\subsection{Potential optimality classification} \label{apA-02}

Potential optimality of a design choice is classified by its Pareto optimality (\cref{th15}). More precisely, the classification is based on the rank $r$ of the point (\cref{th06}). The specification of the Lipschitz constant is tightened from $K_\text{L}>0$ to $K_\text{L}=\varepsilon$ where $\varepsilon$ is a very small positive number, giving
\begin{equation}
	r_j - K_\text{L}d_j \leq r_j - K_\text{L}d_i, \; \forall i\in P, \; K_\text{L}=\varepsilon\label{eq06}.
\end{equation}
This criterion helps to quickly localize regions of optimal solutions, ultimately allowing for a more efficient convergence given in \cref{cn03}.

\begin{figure}[t]
\centering
\includegraphics[scale=1,trim=0 5 0 3]{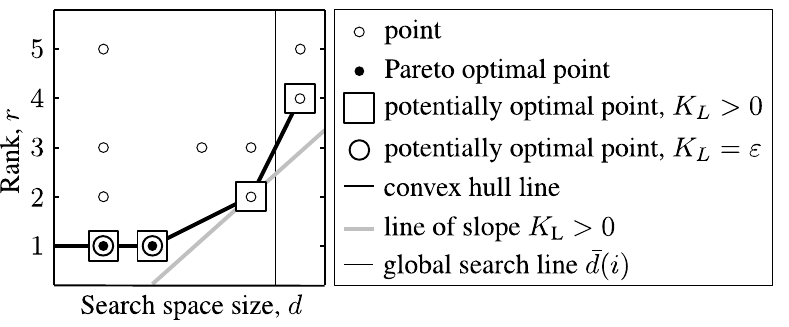}
\caption[]{Potential optimality selection.}
\label{fig01}
\end{figure}

Potential optimality selection can be intuitively illustrated on an $f$-$d$ plot (Fig.~\ref{fig01}). A point $j$ satisfies \eqref{eq05} if there is a line intersecting the point with a gradient $K_\text{L}>0$ such that all other points lie above the line in $f$-$d$ coordinates. Consequently, all potentially optimal points lie on the lower right edge of the convex hull of the points. The tightened requirement to $K_\text{L}=\varepsilon$ in \eqref{eq06} effectively means that only points with the lowest objective values (rank) are chosen to be potentially optimal.

In addition to \eqref{eq06}, the criterion
\begin{equation}
\vspace{-1mm}
	d_j > \overline{d}(i) \label{eq12}
\vspace{-1mm}
\end{equation}
is used to guarantee convergence in \cref{cn01}. $\overline{d}(i)$ is monotonically decreasing with $i$ and $\lim_{i\rightarrow\infty}\overline{d}(i)=0$. This asserts that search spaces that are relatively much bigger are divided, effectively acting as global search. Consequently, eventually all the search space divisions will be divided regardless of satisfaction of~\eqref{eq06}.

\subsection{Search space normalization, bounds, division and size} \label{apA-03}

\begin{figure}[t!]
\centering
\includegraphics[scale=1,trim=0 5 0 0]{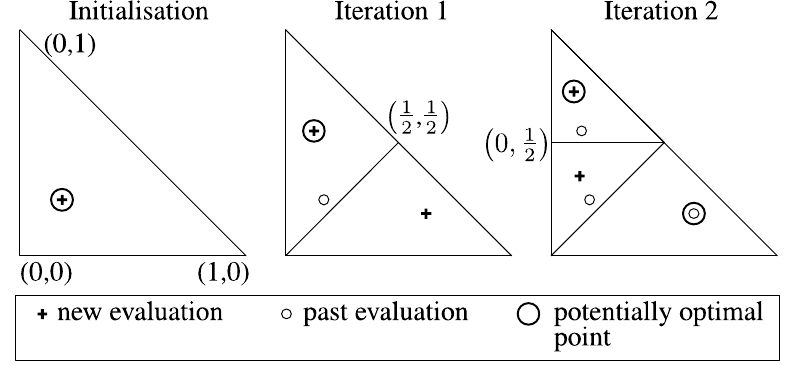}
\caption[]{Simplex division over two iterations.}
\label{fig05}
\end{figure}

The search space is bounded rectangularly, specified by $\underline{\mathbf{p}}$ and $\overline{\mathbf{p}}$ containing the lower- and upper-bounds for each design parameters. A normalized point $\mathbf{c}$ is defined as
\begin{align}
\vspace{-1mm}
	\mathbf{c} &= (\mathbf{p}-\underline{\mathbf{p}})\oslash(\overline{\mathbf{p}}-\underline{\mathbf{p}}),
\\
	\text{so that }\; \mathbf{p} &= \underline{\mathbf{p}} + (\overline{\mathbf{p}}-\underline{\mathbf{p}})\otimes\mathbf{c}. \label{eq15}
\vspace{-1mm}
\end{align}
\Cref{th05} allows for improving search efficiency by focusing on the band $\mathcal{P}_\text{b}$ defined in the theorem. To efficiently locate the specified band, the search will be simplex based (triangular) as illustrated in Fig.~\ref{fig05}. In each iteration, every potentially optimal simplex is divided to form two simplexes of equal size. The size measure $d$ used in criteria \eqref{eq06} and \eqref{eq12} is the longest distance from the center to the vertices of the simplex.

\begin{figure}[t!]
\centering
\includegraphics[scale=1,trim=0 5 0 0]{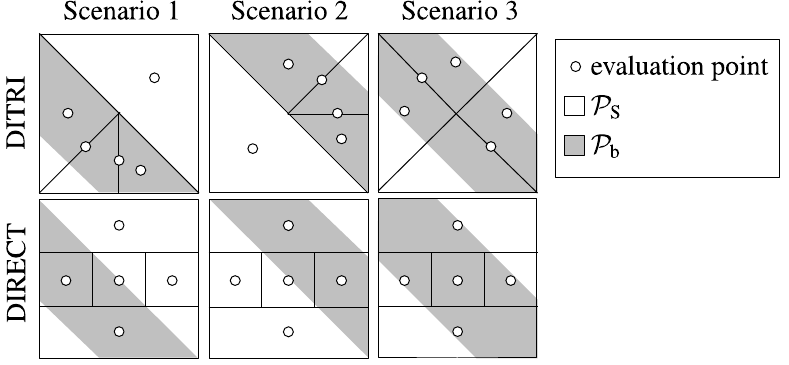}
\caption[]{Comparison between the more adaptable and effective simplex based search compared to rectangular search. Each evaluation point is placed in the middle of the search space. $\mathcal{P}_\text{b}$ is defined in \cref{th05}.}
\label{fig06}
\end{figure}

The efficiency of a simplex based search comes from the fact a simplex is the basic polytope in any $n$-dimension. The approach is contrasted to the classical implementation of Lipschitzian optimization whereby hyper-rectangles are used (DIRECT \cite{Jones:01}). The efficiency of DITRI is demonstrated in Fig.~\ref{fig06}. DITRI does the minimum evaluations (two) per iteration and thus can adjust the search direction more efficiently compared to DIRECT (either two or four). For a given limit on evaluation count, DITRI can more efficiently locate the optimal regions $\mathcal{P}_\text{b}$. This complies with \cref{cn03} for the proposed algorithm.

\vspace{-3mm}
\subsection{Evaluation point location} \label{apA-04}
\vspace{-2mm}

Point evaluation within a simplex is determined stochastically instead of (deterministically) at the center. The uniform  sampling is such that the expected value for the evaluation point chosen is at the center of the simplex,
\vspace{-2mm}
\begin{equation}
	\mathbf{c} = \mathbf{s}_{n_p+1} \label{eq14}
\vspace{-1mm}
\end{equation}
where
$\mathbf{s}_1 \coloneqq b_1$,
$\mathbf{s}_i \coloneqq b_i + \mathsf{U}[0,1]^{1/(i-1)}(\mathbf{s}_{i-1}-b_i),\; \forall i\in\{2,\ldots,n_p+1\}$, and $b_i$ defines the coordinates of the $i^\text{th}$ simplex vertex.


The approach of a random evaluation point allows a faster convergence rate on average than that achieved by a deterministic evaluation. This is possible because the number of instances when a random sample is better placed than the midpoint is on average equal to the instant when it is worse placed. At the case when a random sample is worse placed, its effect would be dominated by the better placed sample and diminished at subsequent iterations.

Finally, integer-valued design parameters are handled simply by shifting the design parameter $p_i$ to the nearest integer, or to the immediate larger integer if it is a half-integer, for any integer-valued design parameter $i$. This fulfills \cref{cn01b}. The number of prediction steps $N$ is the relevant integer-valued design parameter.

The detailed outline of DITRI is given in Algorithm~\ref{al01}.

\begin{algorithm}[t!]
\caption{DITRI}
\begin{algorithmic}[1]
\Require Design bounds $(\underline{\mathbf{p}},\overline{\mathbf{p}})$,\vspace{-.5mm}
\Statex maximum evaluation $\overline{i^\text{ev}}$ and iteration counts $\overline{i^\text{it}}$

\State $i^\text{it} \gets 0$, $i^\text{ev} \gets 0$

\State Let $\mathbf{\Phi}$ contain all evaluated points

\For {the two initial points $i=1$ and $2$}
	\State Set normalized vertices of initial point $\mathbf{b}_i$ 
	\State Evaluate point $\boldsymbol\ell(\mathbf{p}_i)$ in $\mathbf{b}_i$ using \eqref{eq15}, \eqref{eq14}
	\State Record initial point in $\mathbf{\Phi}$
	\State $i^\text{ev} \gets i^\text{ev}+1$
\EndFor

\State $\mathcal{S}_0 \gets \{\mathbf{p}_1,\mathbf{p}_2\}$

\Loop
	\State Find set $O$ of potentially optimal points in $\mathbf{\Phi}$
	\ForAll {$o \in O$}
		\State Divide simplex $\mathbf{b}_o$ to obtain simplices $\mathbf{b}_1^+$, $\mathbf{b}_2^+$
		\For {\textbf{both} $j=1$ and 2}
			\State Evaluate point $\boldsymbol\ell(\mathbf{p}^+)$ in $\mathbf{b}_j^+$ using \eqref{eq15}, \eqref{eq14}
			\State Record new point in $\mathbf{\Phi}$
			\State $i^\text{ev} \gets i^\text{ev}+1$
			\State Update $\mathcal{S}_{i^\text{ev}}$ to be all Pareto points in $\mathbf{\Phi}$
			\State \textbf{if} {$i^\text{ev}\geq\overline{i^\text{ev}}$} \textbf{then} terminate algorithm
		\EndFor
	\EndFor
	\State $i^\text{it} \gets i^\text{it}+1$
	\State \textbf{if} {$i\geq\overline{i^\text{it}}$} \textbf{then} terminate algorithm
\EndLoop
\end{algorithmic}
\label{al01}
\end{algorithm}

\section{Proof of \cref{th18}} \label{apC}

\begin{proof}
Let $P(i^\text{it})$ be the set containing all the evaluated points at iteration number $i^\text{it}$ in Algorithm~\ref{al01}. After some $i_+$ steps ahead, every point in $P(i^\text{it})$ would eventually be classified as potentially optimal either via criterion \eqref{eq06} given its rank, or criterion \eqref{eq12} given that each search space size $d_j$ for all $j\in\{1,\ldots,|P(i^\text{it})|\}$ would become smaller than $\overline{d}(i^\text{it}+i_+)$. Each potential optimal search space will be divided and at least 1 new point will be evaluated after each division so that $P(i^\text{it}+i_+)\supset P(i^\text{it})$. Consequently, as $i_\text{it}\rightarrow\infty$, the algorithm would search the space $\mathcal{P}_\text{s}$ entirely i.e.\ $\lim_{i^\text{it}\rightarrow \infty}P(i^\text{it})=\mathcal{P}_\text{s}$. Any point $\mathbf{p}\in \mathcal{P}_\bullet\subseteq\mathcal{P}_\text{s}$ in the Pareto design set will be evaluated so that $\lim_{i^\text{it}\rightarrow \infty} S_{i^\text{it}}=\mathcal{P}_\bullet$.
\end{proof}

\section*{References}
\bibliographystyle{abbrv}
\bibliography{MODMPC}

\end{document}